\documentclass[12pt]{amsart}

\usepackage{amssymb}
\usepackage{verbatim}

\title[Weakly Closed Unipotent Subgroups]
{Weakly Closed Unipotent Subgroups\\
in Chevalley Groups}

\author[R.\ Guralnick and G.\ R\"{o}hrle]
{Robert Guralnick and Gerhard R\"{o}hrle}

\address
{Department of Mathematics,
University of Southern California,
3520 S. Vermont Ave.
Los Angeles, CA 90089-2532,
USA}
\email{guralnic@usc.edu}

\address
{School of Mathematics, University of Southampton,
Southampton SO17 1BJ, UK}
\email{G.Roehrle@soton.ac.uk}

\dedicatory{Dedicated to Bernd Fischer on the occasion of
his 70th birthday}

\makeatother

\numberwithin{equation}{section}

\newtheorem{thm}[equation]{Theorem}
\newtheorem{lemma}[equation]{Lemma}
\newtheorem{prop}[equation]{Proposition}
\newtheorem{cor}[equation]{Corollary}

\theoremstyle{definition}

\newtheorem{rem}[equation]{Remark}
\newtheorem{ex}[equation]{Example}

\DeclareMathOperator{\Ad}{{Ad}}                               
\DeclareMathOperator{\Char}{{char}}
\DeclareMathOperator{\kernel}{{ker}}
\DeclareMathOperator{\Lie}{{Lie}}                               
\DeclareMathOperator{\rank}{{rank}}                               

\newcommand\fc{\mathfrak c}
\newcommand\fg{\mathfrak g}

\newcommand{\Size}[1]{\left|#1\right|}


\begin{document}

\begin{abstract}
The goal of this note is to classify the 
weakly closed unipotent subgroups in 
the  split Chevalley groups.
In an application we show under some mild assumptions
on the characteristic that 
the Lie algebra of 
a connected simple algebraic group 
fails to be a so called 2F-module.
\end{abstract}

\maketitle

\section{Introduction}
\label{sec:setup}

Let $G$ be a group and let $H \le K$ be subgroups of $G$. The 
subgroup $H$ is said to be \emph{weakly closed in $K$} if $H$ is  
the only $G$-conjugate of itself contained in $K$.
The notion of weak closure has been quite important in finite
group theory.   

The aim of this note is to classify all 
weakly closed unipotent subgroups of a Borel 
subgroup in the split 
Chevalley groups and to obtain partial results for
the finite twisted (or quasi-split) Chevalley groups.
It is well known that unipotent radicals of parabolic subgroups
are weakly closed (Lemma \ref{l:radical}).
Under some mild restrictions on the size of the 
underlying field, we show in our main result that 
the converse holds, i.e., that
a weakly closed unipotent subgroup of a split Chevalley group
is the unipotent radical of a parabolic subgroup 
(Theorem \ref{mainthm}).   
However, even in split Chevalley groups defined over 
very small fields there are other examples (Example \ref{ex1}).

We obtain less complete results for the  twisted groups 
(Theorem \ref{mainthm-twisted}) and show
that there are other examples of weakly closed subgroups no matter
what the field size is in certain cases (Examples \ref{ex3}, \ref{ex4}).
We only indicate examples for the Suzuki and Ree groups.
 
In analogy to the finite group case, 
for $G$ a connected simple algebraic 
group we say that a non-trivial (irreducible)
$G$-module $V$  is a \emph{2F-module} provided 
\begin{equation}
\label{e:2F}
2\dim X + \dim C_V(X) \ge \dim V,
\end{equation}
where $X$ is a (closed) unipotent 
(but not necessarily abelian or connected) 
subgroup of $G$ and 
$C_V(X)$ denotes the subspace of $X$-fixed points of $V$.  

For the concept and 
relevance of 2F-modules in finite group theory,
we refer the reader to 
\cite{as}, \cite{GurMalle}, and \cite{GurMalleII}.
Here the original question is for a finite group $G$ and a given 
absolutely irreducible faithful $G$-module $V$ to find the maximum
of the expression $|X|^2 \cdot |C_V(X)|$, 
when $X$ is a non-trivial elementary abelian unipotent subgroup of $G$,
cf.\ \cite{GurMalle}.

For the finite simple groups there are very few 
2F-modules, see \cite{GurMalle}, \cite{GurMalleII}.
Analogously, we briefly discuss the sparsity 
of 2F-modules for a simple algebraic group $G$ 
(Remark \ref{r:2F}).
In particular, we apply our main theorem to 
show that the adjoint module $\fg = \Lie G$
of $G$ is not a 2F-module
(Corollary \ref{c:2F}). 
This generalizes a result
of Guralnick and Malle, \cite{GurMalle}, \cite{GurMalleII}.

\smallskip

We assume throughout that the groups of Lie type considered
are generated by unipotent elements.

Let $k$ be a field. With the exception of Section 3, 
$G = G(k)$ denotes a split (adjoint) Chevalley group 
in the sense of \cite{Steinberg}.  

Let $T$ be a Cartan subgroup 
of $G$ and $B$ is a Borel subgroup
of $G$ containing
$T$. Let $U \le B$ be the unipotent radical of $B$.
Let $\Psi = \Psi(G,T)$ be the root system of $G$ with respect to $T$
and let $\Pi = \Pi(B)$ be the set of simple roots of $G$
and $\Psi^+ = \Psi(B)$ the set of positive roots of $G$
defined by $B$.
For $\gamma \in \Psi$ we denote the 
root subgroup
defined by $\gamma$ by $U_\gamma$.
For a subgroup $H$ of $G$ we set
$\Psi(H) = \{\beta \in \Psi \mid U_\gamma  \le H\}$.
By $W$ we denote the Weyl group of $G$ with respect to $T$.

Let $P \ge B$ be a parabolic subgroup of $G$.
Then $P$ factors as $P = LP_u$ with some Levi complement $L$
and unipotent radical $P_u$.
In such a decomposition we 
always assume that $L$ is \emph{standard}, 
i.e., that $L$ is generated by $T$ along with the root subgroups of 
a subsystem of $\Psi$ which is generated by a subset of the simple roots,
e.g., see \cite[\S 2.6]{Carter}.

We say that {\em $p$ is a very bad prime for $G$} if $p$ divides one
of the structure constants of the Chevalley commutator relations
for $G$, \cite [p.\ 12]{Steinberg}; that is $2$ (resp.\ $3$) 
is a very bad prime 
for $G$ if $G$ admits a simple factor of type $B_n$, $C_n$, for $n \ge 2$,
$F_4$, or $G_2$ (resp.\ $G_2$);
else there are no very bad primes for $G$.

As general references for Chevalley groups and algebraic groups 
we refer the reader
to \cite{Bo}, \cite{Carter0}, \cite{Carter}, and \cite{Steinberg}.

\section{Weakly Closed Unipotent Subgroups}
\label{sec:weak1}

We maintain the notation and assumptions from the Introduction.
In particular, in this section 
$G = G(k)$ denotes a split Chevalley group.
We first show that the unipotent radicals of parabolic subgroups
are weakly closed.  This is a well known fact; it is 
stated almost in this form in \cite[Lem.\ 4.2]{grodal}.   
This is also proved in \cite[I.2.5]{as}, based on
Alperin--Goldschmidt fusion. 

\begin{lemma}
\label{l:radical}
Let $P \le G $ be a parabolic subgroup of $G$.
Then $P_u$ is weakly closed in $U$.
\end{lemma}

\begin{proof}  
Suppose $P_u^g \le U$ for $g \in G$.
Let $w \in W$
be the minimal length double coset representative
of the $(P, B)$-double coset in $G$ containing $g$. Then $P_u^w \le U$.
Suppose $w \ne 1$. Then
for some simple root $\alpha \in \Pi \setminus \Psi(L)$
the simple reflection $s_\alpha$  is  a prefix of $w$, i.e.,
$w$ has a reduced expression beginning with $s_\alpha$.
Since $U_\beta = U_\alpha ^w \le P_u^w \le U$
and $\beta$ is a  negative root,
this is a contradiction.
Consequently, $w=1$ and thus $g \in P$ and so $P_u^g = P_u$,                   
as desired.
\end{proof}

\begin{lemma}
\label{l:normal}
Let $X \le U$ be weakly closed in $U$.
Set $P = N_G(X)$. 
Then $P$ is a parabolic subgroup of $G$ and $X$ 
is contained in $P_u$.
\end{lemma}

\begin{proof}  
Since $X$ is weakly closed in $U$, we have $B = N_G(U) \le N_G(X) = P$ and 
so $P$ is a parabolic subgroup of $G$.
As $P_u$ is the largest normal unipotent subgroup of $P$, we have
$X \le P_u$.
\end{proof}

If $X$ is a weakly closed subgroup of $U$, 
it is normalized by $T \le B$.
If $k$ is sufficiently large, then 
$X$ is  generated by the root groups 
contained in it, cf.\ \cite{azad}, 
\cite[Prop.\ 14.4.(2a)]{Bo}, 
\cite{Seitz1}, and \cite[Lem.\ 17]{Steinberg}.
In order to ensure that,  
we make the following restrictions on $k$; 
e.g., see \S 5  in \cite{Vavilov1} and in particular 
the references therein.

\begin{itemize}
\item[$(\dag)$] 
$k \ne \mathbb F_2, \mathbb F_4$ in case $G$ is of type $A_2$; \\
$k \ne \mathbb F_2, \mathbb F_3$ in case $G$ is of type $A_3$,
$B_n$, $C_n$, for $n \ge 2$, $D_n$, for $n \ge 3$, $F_4$, or $G_2$; \\ 
$k \ne \mathbb F_2$ in case $G$ is of type $A_n$,
for $n \ge 4$, $E_6$, $E_7$, or $E_8$;\\
$k$ is perfect if $\Char k = 2$ and $G$ is of type $C_n$, for $n \ge 1$.
\end{itemize}

Our main result gives a converse to Lemma \ref{l:radical} assuming 
$(\dag)$.

\begin{thm} 
\label{mainthm}
Assume $(\dag)$. 
Suppose  $X \le U$ is weakly closed in $U$. Set $P = N_G(X)$.
Then $X = P_u$.
\end{thm}

\begin{proof}
By Lemma \ref{l:normal}, $P$ is a parabolic subgroup of $G$
and $X \le P_u$.

Since $X$ is normalized by $T\le B$, 
the restrictions on $k$ in $(\dag)$ ensure that  
$X$ is generated by the root subgroups contained in $X$.
Since $X$ is normalized by $B$, it follows from the
commutator relations of $G$ (cf.\ \cite[p. 30, Lem.\ 33]{Steinberg})
that $\Psi(X)$ is a closed subset of $\Psi^+$, in the sense of 
\cite[p.\ 24]{Steinberg}.
Consequently, we have $X = \prod U_\beta$,
where the product is taken in any fixed order over $\Psi(X)$, 
cf.\ \cite[Lem.\ 17]{Steinberg}.

Now suppose that there is a simple root $\alpha \in \Pi \cap  \Psi(P_u)$
such that $U_\alpha \nleq X$.
Then $U_\alpha \cap X = \{1\}$, since $X$ is generated by root subgroups
and distinct root subgroups intersect trivially; the latter follows from the 
uniqueness of factorization in the product decomposition $X = \prod U_\beta$,
\cite[Lem.\ 17]{Steinberg}.
Then $X^{s_\alpha} \le U$, since $s_\alpha$ permutes
$\Psi^+\setminus\{\alpha\}$.
Thus, as $X$ is weakly closed in $P_u$,
we have $s_\alpha \in N_G(X) = P$.
But for a simple root $\alpha$, we have 
$s_\alpha \in P$ if and only if $\alpha \in \Psi(L)$ 
if and only if $\alpha \not\in \Psi(P_u)$,
a contradiction.
Consequently,
$\Pi \cap  \Psi(P_u) \subseteq \Psi(X)$.
Since $P_u$ is the normal closure in $P$ of the root groups relative to
the simple roots in $\Psi(P_u)$,
e.g., see \cite[Prop.\ 2.10, Rk.\ 2.13]{ParkerRoehrle},
we derive that $P_u \le X$.
\end{proof}

We recall a well-known fact concerning 
regular unipotent elements.

\begin{rem}
\label{r:regular}
Let $G$ be a reductive algebraic group. 
A unipotent element $u$ of $G$ is called \emph{regular}
provided $\dim C_G(u)$ is minimal possible among unipotent elements
in $G$. A regular unipotent element is contained in a 
unique Borel subgroup of $G$, see \cite[Prop.\ 5.1.3]{Carter}.
Let $F$ be a Frobenius endomorphism of $G$ so that the subgroup of 
$F$-fixed points $G^F$ of $G$ is a finite group of 
Lie type.
Let $u \in G^F$ be regular unipotent in $G$.
Since any Borel subgroup of $G^F$ is the 
fixed point subgroup of a unique $F$-stable Borel subgroup of $G$,
the uniqueness result for $G$ just quoted implies that 
$u$ is in a unique Borel subgroup of $G^F$.
\end{rem}

The following example shows that the hypothesis $(\dag)$ of 
Theorem \ref{mainthm} is necessary.

\begin{ex}
\label{ex1}  
Let $G$ be a split simple Chevalley group
over the field of $2$ elements of rank at least $2$.   
Note that $U=B$. Let $Y$ be the subgroup
of $U$ generated by the root subgroups $U_\gamma$ relative to all
the non-simple positive roots, i.e.,  $\gamma \in \Psi^+\setminus \Pi$.  
So every subgroup of $G$ 
between $U$ and $Y$ is normal in $B$.  
Let $u \in U$ be a regular unipotent element.
Note that this determines the coset $uY$ uniquely.
Let $X$ be the subgroup of $G$ generated by $u$ and $Y$.   
We claim that $X$ is weakly closed in $U$.   
Since $u$ is regular unipotent, it is contained in no
other Borel subgroup of $G$, cf.\ Remark \ref{r:regular}, 
and so the same is true for $X$. 
So if $X^g \le U$, then $X \le U^{g^{-1}}$ and so $g \in N_G(U)=U$, whence 
$X^g=X$. Thus $X$ is weakly closed in $U$.
Since $\rank G > 1$, it follows that $X$ is not the unipotent
radical of any parabolic subgroup of $G$.
\end{ex}

One can construct in a similar way examples for all cases of 
split groups when $(\dag)$ fails.

\section{Weakly Closed Subgroups in Finite Twisted Groups}

We note that Lemma \ref{l:radical} also holds for the finite
twisted Chevalley groups and the proof goes through verbatim
only involving the $(B,N)$-pair structure of the underlying group.
We record this:

\begin{lemma}  
\label{2:radical-twisted case}
Let $G$ be a finite simple Chevalley group and 
$P \le G $ be a parabolic subgroup of $G$.
Then $P_u$ is weakly closed in $U$.
\end{lemma}

This is also proved in \cite{as} and \cite{grodal}.  We sketch 
a proof of Lemma \ref{2:radical-twisted case} 
for classical groups in Lemma \ref{3radical-classical case} below
that is quite different from
the other proofs mentioned.  For groups of rank $1$, there is nothing
to prove.  We do not complete the argument for the 
exceptional groups, 
but the proof 
of Lemma \ref{3radical-classical case} below does show 
that it suffices to check the statement for unipotent radicals 
of the maximal parabolic subgroups.
By a classical group, we mean a linear, unitary, symplectic or orthogonal
group.

We first recall some general properties of weakly closed subgroups
for finite groups (with obvious modifications 
Lemma \ref{general weakly closed} also applies to unipotent
subgroups of algebraic groups).

\begin{lemma}
\label{general weakly closed}
Let $G$ be a finite group with $U$ a Sylow $p$-subgroup of $G$.
Let $B=N_G(U)$.   Let $X$ be a normal subgroup of $U$.
The following are equivalent:
\begin{enumerate}
\item $X$ is weakly closed in $B$;
\item $X$ is weakly closed in $P:=N_G(X)$;
\item $X$ has a unique fixed point on $G/P$.
\end{enumerate}
Moreover, if any of these conditions holds, then $P=N_G(P)$.
\end{lemma}

\begin{proof}   
Assume that $X^g \le P$ for some $g \in G$.  
Then the subgroup of $G$ generated by $X$ and $X^g$
is a $p$-subgroup of $P$ and so by conjugating, we may assume that
$X^g \le U$.  Thus, (1) implies (2).  Since $U \le P$, 
(2) implies (1).


Note that $X^g \le P$ if and only if $X$ fixes the point 
$gP$ in $G/P$.
Thus, if $X$ is weakly closed in $P$ and if 
$X$ fixes $gP$ in $G/P$, then $X=X^g$, 
and so $g \in N_G(X) = P$.
Thus (3) follows from (2).


Now assume that (3) holds.  
Since $X$ fixes the point $P$ in $G/P$, 
if $X^g \le P$, we have $gP = P$ and so $g \in P$.
Thus (2) holds.  

Finally, if $g$ normalizes $P$, then $X^g \le P$, 
whence $g$ normalizes $X$ and so is in $P$.
So, the last assertion follows.
\end{proof} 

\begin{lemma}
\label{3radical-classical case}
Let $G$ be a simple classical group over a finite field.
If $P$ is a parabolic subgroup of $G$,
then $P_u$ is weakly closed in $B$.
\end{lemma}

\begin{proof} 
We argue by induction on the rank of $G$.  If $\rank G =1$, 
then $U$ is the unique Sylow $p$-subgroup of $B$, 
and hence is weakly closed in $B$.

Let $N$ be the natural module of $G$.
First consider the case that $P$ is a maximal parabolic subgroup of $G$.
Then $P$ is the stabilizer of a totally singular $m$-subspace $M$ 
of $N$ for
some $m \le \dim N/2$ and $P_u$ is the subgroup of $G$ acting trivially
on both $M$ and $M^{\perp}/M$ (for the linear case we have $M^{\perp} =  N$).
Note that $M$ is precisely the set of fixed singular
vectors for $P_u$.

We claim that $M$ is the unique $P_u$-invariant totally singular
subspace of dimension $m$ (of the given type in the case of 
orthogonal groups).  We show this by induction on $m$.
The case $m = 1$ is clear. Let $m > 1$.
Now let $V$ be a $P_u$-invariant totally singular
subspace of $N$ of dimension $m$.  
Then $P_u$ fixes some non-zero vector 
$v \in V$ and so $v \in M$.  So $V$ is contained in 
$\langle v\rangle^{\perp}$.
By induction, $P_u$ fixes a unique totally isotropic 
subspace of dimension
$m-1$ (of the given type) in $\langle v\rangle^{\perp}/\langle v\rangle$,
whence the claim.

Thus,  $P_u$ has a unique fixed point on $G/P$, whence $P_u$ is 
weakly closed in $B$, thanks to Lemma \ref{general weakly closed}.

Suppose that $P$ is not maximal.  Let $Q$ be a maximal parabolic
subgroup of $G$ containing $P$.  
Then $P_u/Q_u$ is the unipotent radical of $P/Q_u$
in $Q/Q_u$ and the latter is a central
product of a torus and some number of smaller classical groups. 
Since $P/P_u$ contains the torus and is a central product of 
parabolic subgroups
in each factor, it follows by induction that $P_u/Q_u$ is weakly closed in
$B/Q_u$.   

Observe that $P_u^g \le B$ implies that $Q_u^g \le B$.
Since $Q_u$ is weakly closed in $B$, by the case above,
it follows that $g \in Q$.  
Thus, as $P_u^g \le B$ implies that $g \in Q$
and since $P_u/Q_u$ is weakly closed in
$B/Q_u$, we have $(P_u/Q_u)^{gQ_u} = P_u/Q_u$ and thus
we obtain $P_u^g=P_u$, as desired.
\end{proof}

Our next examples show that in the twisted groups there
are additional instances of weakly closed unipotent subgroups
for \emph{all} finite fields.  The proof of
Lemma \ref{3radical-classical case}
shows that the critical case is that of a maximal parabolic subgroup.

\begin{ex}  
\label{ex2}
Let $G = \mathrm{U}_3(q)$.  Take $X=Z(U)$.  So $|U|=q^3$ and
$|X|=q$.  Note that $X$ consists of all transvections in $U$ and so clearly
it is weakly closed in $U$.
\end{ex}

\begin{rem}
\label{r:ree}
Note that a minor variation of the previous example shows
that for all the twisted rank $1$ groups there are proper
weakly closed subgroups of $U$.  Similarly, this holds for
Ree groups of type $F_4$.  
\end{rem}

Since $\mathrm{U}_3(q)$ is a Levi subgroup of a 
parabolic subgroup of $\mathrm{U}_{2m+1}(q)$ for
all $m > 1$, we can use Example \ref{ex2} to give other examples.

\begin{ex} 
\label{ex3}
Let $G=\mathrm{U}_{2m+1}(q)$ with $m > 1$.  Let $P$ be a minimal
parabolic subgroup of $G$ with unipotent radical $R$ such
that the derived subgroup of $P/R$ is a $3$-dimensional unitary group.
Let $X$ be the subgroup of $P$ such that $X/R$ is the center of 
$U/R$.  Since $R$ is weakly closed in $U$,
by Lemma \ref{2:radical-twisted case}, 
it follows that $X^g \le U$
implies that $R^g \le U$ and so $R^g=R$, i.e.\ $g \in P$.  By 
Example \ref{ex2}, this implies that $X^g=X$, as required.
\end{ex}

In fact, there are additional examples of weakly closed subgroups $X$
where $N_G(X)$ is even a maximal parabolic subgroup of $G$.

\begin{ex}
\label{ex4} 
Let $G=\mathrm{U}_{2m+1}(q)$ with $m \ge 1$.  Let $P \ge U$ 
be the parabolic subgroup of $G$
that is the stabilizer of a totally singular $m$-subspace of 
the natural module $N$ of $G$.  Let $X$ be
the derived subgroup of $P_u$.  
We claim that  $X$ is weakly closed in $U$
and proper in $P_u$.

Since $P_u$ is nilpotent, $X$ is proper in $P_u$.  Let $V$ be the
set of fixed points of $P$ on $N$.
Then $V$ is $m$-dimensional and totally singular.

We claim that $V$ is the only totally singular $m$-space 
left invariant by $X$.  We prove this by induction on $m$.
The case $m=1$ is clear.  Let $m > 1$. Note that the subspace of fixed points
of $X$  on $N$ is $V^{\perp}$ and has dimension $m+1$.  Moreover every vector
in $V^{\perp} \setminus{V}$ is non-singular.  So if $X$
leaves invariant a totally singular $m$-space $V'$, then $X$ is trivial
on some $1$-space $V_1 \le V' \cap V_1^{\perp}$. Now view $X$ acting on
$V_1^{\perp}/V_1$.  
This gives a homomorphism of $X$ into $\mathrm{U}_{2m-1}(q)$.
The image of $X$ is precisely the derived subgroup of the unipotent radical
of the maximal parabolic subgroup stabilizing $V/V_1$, whence by induction
$V'/V_1 = V/V_1$ and so $V=V'$, as claimed.

So $X$ has a unique fixed point on $G/P$, whence $X$ is weakly closed
in $P$ and thus also in $B$, by Lemma \ref{general weakly closed}.
\end{ex}

\begin{rem}
One might wonder why the previous example does not extend to odd
dimensional orthogonal groups (or twisted orthogonal groups).  
The problem is that the first step in the induction fails.   
\end{rem}

Examples \ref{ex3} and \ref{ex4} show that there are additional 
weakly closed unipotent subgroups 
in the odd-dimensional unitary groups (over any finite field).
However, our main result of this section shows 
that for the remaining families of twisted groups, 
we have the same result as Theorem \ref{mainthm},
as long as the field size is sufficiently large.


\begin{thm}
\label{mainthm-twisted}   
Assume that $k$ is a finite field with $|k| > 5$.
Let $G={^2}A_{2m+1}(k)$, for $m \ge 1$, 
${^2}D_n(k)$ for $n \ge 4$, 
${^3}D_4(k)$,  or ${^2}E_6(k)$.
Let $U$ be the unipotent radical of a Borel subgroup $B$ of $G$.  If
$X \le U$ is weakly closed in $B$, then $X=P_u$ for some parabolic subgroup
$P$ of $G$.
\end{thm}

\begin{proof}  Let $T$ be a maximal torus contained in $B$.  So $X$
is normalized by $T$.
 The assumption on $|k|$ guarantees that $X$ is a product of root subgroups,
 see \cite[Lem.\ 3]{Seitz1}.  
 Moreover, in all the cases considered, the root subgroups
 are abelian (and can be identified with $k$ or the quadratic extension of $k$
 or, in the case of ${^3}D_4(k)$, the cubic extension of $k$) 
 and the intersection of
 any two distinct root subgroups is trivial.   Let $P=N_G(X)$, a parabolic
 subgroup of $G$. So $X \le P_u$.  As in the previous section, in this case, we
 see that $X$ must contain all the root subgroups corresponding to simple
 roots in $\Psi(P_u)$ 
 (this uses the fact that root subgroups relative to distinct roots
 intersect
 trivially and that the simple reflections $s_\alpha$ in 
 the Weyl group of $G$ preserve
 the set of positive roots other than the simple root $\alpha$).  
 As in the split case, we see that the normal closure of the simple root
 subgroups contained in $P_u$ is all of $P_u$.  Thus, $X=P_u$.
\end{proof}

\section{Centralizers of Weakly Closed Unipotent Groups}
\label{sec:wcs}

For the remainder of the note we 
assume that $G$ is a connected, simple algebraic group
and that $k$ is algebraically closed.
For a (closed) subgroup $H$ of $G$ we denote the
identity component of $H$ by $H^0$ and note that 
$\dim H = \dim H^0$ meaning dimension as an algebraic variety.

In this section we show that the adjoint module $\fg$ of $G$
is not a 2F-module, in the sense of \eqref{e:2F}.
This extends a result due to Guralnick and Malle \cite{GurMalle}.
Moreover, we show that there are very few possibilities for
a 2F-module for dimension reasons.

We only prove the results for algebraic groups.  We leave
it to the reader to prove the same statements for the split
Chevalley groups satisfying $(\dag)$.

\begin{lemma}
\label{l2}
Let  $P$ be a proper parabolic subgroup of $G$.
Then $C_G(P_u) \le P_uZ(G)$. So
$C_G(P_u)^0 \le P_u$ and $C_G(P_u)^0 = Z(P_u)$.
\end{lemma}           
                    
\begin{proof}
Let $L \ge T$ be the standard Levi complement of $P$.
We first show that $C_G(P_u)^0$ contains
no semisimple elements. Suppose that there is a non-trivial
torus $S$ centralizing $P_u$.  Then $S$ is conjugate to
a subtorus of $T$ in $P$ and since $C_G(P_u)$ is normal in $P$,
there is no loss in assuming that $S \le T$.  Then
$S$ centralizes each root subgroup $U_{\alpha} \le P_u$ and
so also $U_{-\alpha}$.  However, as $G$ is simple, it is generated by
$P_u$ and the unipotent radical of the parabolic subgroup of $G$ 
opposite to $P$ (with respect to $L$), 
e.g.\ see \cite[Prop.\ 4.11]{boreltits}.  
Thus, $S \le Z(G)$; a contradiction.
It follows that $C_G(P_u)^0$ is a normal unipotent subgroup of $P$ and so
contained in $P_u$.  Thus, $C_G(P_u)^0=Z(P_u)$.  
The only finite normal subgroups of $L$ are contained in $Z(L) \le T$
and arguing as above, we see that $C_T(P_u)=Z(G)$.
Thus we have $C_G(P_u) \le B$.
Moreover, since $C_G(P_u)$ is normal in $P$, 
we have $C_G(P_u) \le P_u Z(L)$.
By the action of $T$ on the root subgroups of $U$ and
by the 
commutator relations we see that
\[
C_G(P_u) = C_G(P_u) \cap P_u Z(L) = (C_G(P_u) \cap P_u)(C_G(P_u) \cap Z(L)).
\]
Note that $C_G(P_u) \cap Z(L) = \cap_{\alpha \in \Pi} \kernel \alpha = Z(G)$
and since $G$ is simple, the latter is finite. So $C_G(P_u)\le P_u Z(G)$, and
$C_G(P_u)^0 \le P_u$, as claimed.
\end{proof}

\begin{rem}
\label{rem:rich-bound}
Let  $P$ be a parabolic subgroup of $G$.
There is a natural bound for $\dim C_G(P_u)$ 
stemming from Richardson's Dense Orbit
Theorem, e.g., see \cite[\S 5.2]{Carter}.
There is a conjugacy class $C$ of $P$ in $P_u$ which is open dense in $P_u$.
It turns out that for any $x$ in $C$ we get
$C_G(x)^0 = C_P(x)^0$, cf.\ \cite[Cor.\ 5.2.2]{Carter}, and thus 
$\dim C_G(x) = \dim C_P(x)$.
For any $x \in C$ we clearly have 
$C_G(P_u) \le C_G(x)$ and thus, since $\dim C + \dim C_G(x) = \dim P$
and $\dim C = \dim P_u$, we obtain
\[ \dim P_u + \dim C_G(P_u) \le \dim P.\]
The existence of such a dense $P$-orbit is part of a 
fundamental theorem due to R.W.\  Richardson \cite{richardson2}.
The proof relies on the fact that the number of 
unipotent classes of $G$ is finite. This was first proved also by 
Richardson  under some mild restrictions on the characteristic of the 
ground field \cite{richardson1}; these were removed subsequently 
by G.\ Lusztig in \cite{lusztig}. 
\end{rem}

We can improve the bound from Remark \ref{rem:rich-bound}
as follows.

\begin{prop}
\label{p2}
Let $P$ be a proper parabolic subgroup of $G$. Then
\[ \dim P_u + \dim C_G(P_u) \le \dim B.\]
\end{prop}

\begin{proof}
By Lemma  \ref{l2}, 
$C_G(P_u)^0$ is a connected normal abelian subgroup of $U$.
Thus, by \cite[Thm.\ 1.1]{Roehrle},
there are only finitely many $B$-orbits on  $C_G(P_u)^0$,
and consequently there is a dense such orbit.
Thus $\dim C_G(P_u) + \dim C_B(x) = \dim B$ for some $x \in C_G(P_u)^0$.
Finally, since $P_u \le  C_B(x)$, 
the desired inequality follows.
\end{proof}

If $H$ is a (closed) subgroup
of $G$, we define 
\[ f(H) = 2 \dim H + \dim C_G(H).\]
Our next result is a restatement of Lemma 2.1 from \cite{GurMalle}
in our context and the proof is identical, see also
\cite{chermakdelgado}.

\begin{lemma}
\label{2.1} 
Let $H \le M \le G$ and $K \le G$ be subgroups of $G$. 
Suppose that $f(H)$ is maximal among all the subgroups of $M$. 
Let $D$ be the (algebraic) subgroup of $G$ generated by $H$ and $K$.
Then $f(D) \ge f(K)$.
\end{lemma}

\begin{proof} 
Note that
$ \dim D \ge \dim HK = \dim H + \dim K - \dim (H \cap K)$.
Here $HK$ need not be a subgroup of $G$.  Since
$C_G(H \cap K)$ contains both $C_G(H)$ and $C_G(K)$, 
and $C_G(H) \cap C_G(K) = C_G(D)$,
we similarly obtain 
$\dim C_G(D) \ge \dim C_G(H) + \dim C_G(K) - \dim C_G(H \cap K)$.
Thus, we have
 \begin{align*}
f(D) & = 2 \dim D + \dim C_G(D) \\
   & \ge 2(\dim H + \dim K)+\dim C_G(H) + \dim C_G(K) \\
& \qquad\quad - 2 \dim (H \cap K) - \dim C_G(H \cap K)\\
    & = f(H) + f(K) - f(H \cap K) 
 \end{align*}
and since $f(H) \ge f(H \cap K)$,
by maximality of $f(H)$ among the subgroups of $M$, 
the lemma follows.
\end{proof}

For $X$ a subgroup of $U$ let $\widehat X$ denote the 
\emph{weak closure of $X$ in $U$}, that is the smallest
weakly closed  subgroup of $U$ containing $X$ (so 
$\widehat X$ is the subgroup of $U$ generated by all conjugates
of $X$ contained in $U$).

Note that since weakly closed unipotent subgroups are unipotent radicals
of parabolic subgroups of $G$, they in fact are closed.  So we could define
the weak closure of $X$ to be the algebraic group generated by
the conjugates of $X$ contained in $U$.

\begin{cor}
\label{c:weak}
Let $X \le U$ with $f(X)$ maximal 
among all subgroups
of $U$, then $f(X)=f(\widehat X)$, 
where $\widehat X$ is the weak closure of $X$ in $U$.
\end{cor}

\begin{proof}   
Note  that there are
finitely many conjugates $X_i$, $1 \le i \le m$, of $X$ so that
the group generated by $X_1, \ldots, X_m$ has the same centralizer
as $\widehat{X}$.   
Let $Y_i$ be the 
group generated by $X_1, \ldots, X_i$. We show that $f(Y_i)=f(X)$
for all $i$ and this proves the corollary.

This is clear for $i=1$, as $f$ is constant on conjugates.  
Let $i > 1$.
Inductively, we have $f(Y_{i-1}) = f(X)$ and since this is maximal
among all the subgroups of $U$,
we can apply the previous lemma with 
$H = Y_{i-1}$, $M = U$, and $K = X_i$ to conclude that
$f(Y_i) \ge f(X)$ which, by maximality, gives $f(X) = f(Y_i)$.
Thus, $f(\widehat{X})=f(Y_m)=f(X)$.
\end{proof}

\begin{thm} 
\label{2F}  
Let $G$ be a simple algebraic group of rank $r \ge 2$.
In case $\Char k$ is a very bad prime for $G$ assume that 
$r \ge 3$.  Let $X$ be a non-trivial
unipotent subgroup of $G$.  Then 
\[2 \dim X + \dim C_G(X) < \dim G.\]
\end{thm}


\begin{proof}  
Choose $X \le U$ with $f(X)$ maximal.  Choose $X$ to be non-trivial
if possible.  If this is not possible, then $f(X) < f(\{1\}) = \dim G$
and the result follows.

It follows from Corollary \ref{c:weak} that 
the maximum is achieved on the weak closure of $X$ in  $U$.
Thus, thanks to Theorem \ref{mainthm}, we may choose
a parabolic subgroup $P$ of $G$ with $X \le P_u$ and $f(X)=f(P_u)$.
Since $P$ is proper, 
it follows from Proposition \ref{p2} that 
\[2 \dim P_u + \dim C_G(P_u) \le \dim P_u + \dim B < \dim G,\]
unless $P = B$.
If $\Char p$ is not a very bad prime of $G$, then
$\dim Z(U) = 1$, else $\dim Z(U) = 2$.
Thus, by the hypotheses on $r$, we have
\[2 \dim U + \dim Z(U) = \Size{\Psi} + \dim Z(U) < 
\Size{\Psi} + r = \dim G,\]
giving the desired strict 
inequality also for $P = B$. 
\end{proof} 

We next prove an analogue of Theorem \ref{2F} 
for $\fg = \Lie G$, the Lie algebra of $G$; this says 
that $\fg$ is not a 2F-module of $G$, cf.\ \eqref{e:2F}.

For $H$ a (closed) subgroup of $G$, let 
$\fc_{\fg}(H) := \{y \in \fg \mid \Ad(h) y = y \ \text{for all} \  h \in H\}$
denote the subspace of $\fg$ of $\Ad(H)$-fixed points of $\fg$.

We define a function similar  to $f$ above: 
for $H$ a  (closed) subgroup of $G$ set 
\[
f_\fg(H) := 2\dim H + \dim \fc_\fg(H),  
\]
where we use the centralizer in the Lie algebra instead 
of the group.
In general, we have $\Lie C_G(H) \le \fc_\fg(H)$ and thus
$f(H) \le f_\fg(H)$ for any subgroup $H$ of $G$. We have equality 
precisely when the scheme-theoretic centralizer of $H$ in $G$ is smooth.

Using dimension arguments of subalgebras instead of subgroups, 
one readily checks that the proof of Lemma \ref{2.1} also 
applies for $f_\fg$ in place of $f$ (with essentially no change).
Further, since $\fc_\fg(H^0) \ge \fc_\fg(H)$, the proof of 
Corollary \ref{c:weak} is also valid for 
$f_\fg$ in place of $f$.

\begin{cor} 
\label{c:2F}  
Let $G$ be a simple algebraic group of rank $r \ge 2$.
Let $X$ be a non-trivial
unipotent subgroup of $G$. Then 
\[2 \dim X + \dim \fc_{\fg}(X) < \dim \fg.\]
\end{cor}


\begin{proof}
We argue as in the proof of Theorem \ref{2F}.
Choose $X$ with $f_\fg(X)$ maximal with $X$ non-trivial 
(if this is not possible, then 
$f_\fg(X) < f_\fg(\{1\}) =\dim \fg$ for every non-trivial
$X$ and the result holds).  
It follows from the $f_\fg$-analogue of Corollary \ref{c:weak} that 
the maximum is achieved on a (non-trivial)  weakly closed subgroup of $U$.
By Theorem \ref{mainthm}, we have
$X = P_u$ for some proper parabolic subgroup $P$ of $G$.

We show that under our assumptions 
$\Lie C_G(P_u) = \fc_\fg(P_u)$, i.e., that
the scheme-theoretic centralizer of $P_u$ in $G$ is smooth.
In particular, we then obtain that 
$\dim  C_G(P_u) = \dim \fc_\fg(P_u)$ and the desired 
result is immediate from Theorem \ref{2F}.

Since $P \ne G$, 
it  follows from Lemma \ref{l2} that $C_G(P_u)^0 \le P_u \le U$.
And since $P_u$ is $T$-stable, so is $C_G(P_u)^0$.
It thus follows from \cite[Prop.\  14.4(2a)]{Bo} that 
$C_G(P_u)^0 = \prod U_\gamma$, where
the product is taken over the set of roots
$\Gamma := \Psi(C_G(P_u))$.
Thanks to the commutator relations and our restrictions on the 
characteristic, 
$\Gamma =\{\gamma \in \Psi(P_u) \mid \gamma + \beta \not \in \Psi^+\ \forall\ 
\beta\in \Psi(P_u)\}$. This is a closed subset of $\Psi^+$, in the sense of 
\cite[p.\ 24]{Steinberg}.
Likewise, $\fc_\fg(P_u)$ is $\Ad(T)$-stable and thus, 
$\fc_\fg(P_u)$ is a sum of root spaces in $\fg$,
cf.\ \cite[Prop.\  13.20]{Bo}.
Because of the restrictions on $\Char k$, 
there are no degeneracies in the adjoint action of
root elements on root spaces of $\fg$, cf.\ \cite[\S 4.3]{Carter0}.
Consequently, $\fc_\fg(P_u) = \oplus \fg_\gamma$, where the sum is over the 
same set $\Gamma$ defined above.
in particular, we have 
$\Lie C_G(P_u) = \fc_\fg(P_u)$, as claimed.
\end{proof}





\begin{rem}
\label{r:rank1}
If $\rank G = 1$, then the inequalities in Theorem \ref{2F} and
Corollary \ref{c:2F} are clearly 
still valid provided $X$ is a non-trivial finite unipotent 
subgroup; else, of course, we get equality.  
\end{rem}

\begin{rem}
\label{r:2F}
We can consider the same question for any rational $G$-module $V$.
Assume that $V$ is irreducible.  For any (closed) subgroup $H$ of $G$, define
\[
f_V(H):= 2\dim H + \dim C_V(H). 
\]
The question is when there exists
a non-trivial unipotent subgroup $X$ of $G$ with $f_V(X) \ge \dim V$.
As for the Lie algebra case, 
it is straightforward to check that the $f_V$-analogues of 
Lemma \ref{2.1} and Corollary \ref{c:weak} also hold 
with essentially identical proofs.
This then shows that $f_V(X) \le f_V(G) = 2\dim G$, 
since $V$ is irreducible.
So if $V$ is a 2F-module for $G$, i.e., if 
$f_V(X) \ge \dim V$ (cf.\ \eqref{e:2F}), 
then necessarily $\dim V \le 2 \dim G$.
On the other hand, since $V$ is irreducible, we obtain
$f_V(U) = 2\dim U + 1= \dim G - r + 1$.   
So the existence
of such an $X$ is only open for the case 
\[ \dim G - r + 1 < \dim V \le 2\dim G.\]
There are very few irreducible $G$-modules with dimension 
in this range (see \cite{GurMalleII}).  
We have dealt with the adjoint module above.   By the weak
closure result, i.e., 
the $f_V$-analogue of Corollary \ref{c:weak}
and Theorem \ref{mainthm}, we just need to compute 
$f_V(P_u)$ for each parabolic subgroup $P$ of $G$
for the few remaining cases for $V$. We leave the details to
the reader.
\end{rem}



\bigskip 
{\bf Acknowledgments}:
We thank  M.\ Aschbacher, 
G.\ Seitz, S.\ Smith and N.\ Vavilov for helpful comments.
The first author was partially supported by the NSF grant DMS 0140578.
Part of this paper was written during a visit of the first author
to the Institute for Advanced Study (Princeton) and a visit of the second 
author to the Max-Planck-Institute for Mathematics (Bonn).

\end{document}